# ON REPRESENTATIONS OF THE LIE SUPERALGEBRA $P(n)$

VERA SERGANOVA

ABSTRACT. We introduce a new way to study representations of the Lie superalgebra $p(n)$. Since the center of the universal enveloping algebra $U$ acts trivially on all irreducible representations, we suggest to study the quotient algebra $\bar{U}$ by the radical of $U$. We show that $\bar{U}$ has a large center which separates typical finite-dimensional irreducible representations. We give a description of $\bar{U}$ factored by a generic central character. Using this description we obtain character formulae of generic (infinite-dimensional) irreducible representations. We also describe some geometric properties of the supervariety $\operatorname{Spec} \operatorname{Gr} \bar{U}$ in the coadjoint representation.

## 1. INTRODUCTION

In this paper we study representations of the classical Lie superalgebra $p(n)$ introduced in [3]. Let us recall the definitions. Denote by $\widetilde{p}(n)$ the Lie superalgebra of all endomorphisms of $(n|n)$-dimensional vector superspace which preserve a fixed odd non-degenerate symmetric form. The Lie superalgebra $p(n)$ is the commutator of $\widetilde{p}(n)$. It is simple when $n \geq 2$.

The main difficulty in representation theory of $\widetilde{p}(n)$ and $p(n)$ is that the center of the universal enveloping superalgebra of $\widetilde{p}(n)$ is trivial. The center of the universal enveloping of $p(n)$ is non-trivial see [1], but it acts by the same central character on all irreducible representations. Therefore the powerful technique of central characters seems not applicable. We show that it still can work after some reduction.

The universal enveloping algebra $U$ of $\widetilde{p}(n)$ has a non-trivial radical $I$ as was shown in [4]. Let $\bar{U}$ be the quotient algebra $U/I$. We suggest to study the representations of $\bar{U}$ instead of $U$. Using this approach we do not lose any information about irreducible representations but in fact gain some additional information. Our main result is that $\bar{U}$ has a non-trivial center $Z$, which separates typical finite-dimensional irreducible representations. In Theorem 4.8 below we describe $Z$ precisely. Furthermore we show that for a generic central character $\chi : Z \to \mathbb{C}$ the algebra $\bar{U}/\bar{U} \operatorname{Ker} \chi$ is isomorphic to the matrix algebra over $U(\mathfrak{g}_0)/U(\mathfrak{g}_0) \operatorname{Ker} \chi_0$ for a certain central character $\chi_0$ of $U(\mathfrak{g}_0)$. This allows to write character formulae for a generic irreducible weight representation.

The same method is applicable to $p(n)$.

In the last section we suggest some geometric motivation of our construction. We realize the graded algebra $\operatorname{Gr} \bar{U}$ as the algebra of functions on a Poisson subvariety

---

$\mathcal{X}$ in the coadjoint representation $\mathfrak{g}^*$. Let $Q$ be the closure of the union of coadjoint orbits of all even elements. It is a specific feature of supergeometry that $Q$ does not coincide with $\mathfrak{g}^*$ but is a closed subvariety with singularities. We believe that $\mathcal{X}$ coincides with $Q$, however we are able to prove only that $Q$ is an irreducible component of $\mathcal{X}$.

I would like to thank M. Gorelik and I. Musson for very useful discussions and remarks. I am also thankful to M. Duflo and I. Zakharevich for corrections to preliminary version of the paper.

## 2. Preliminaries

Let $\mathfrak{g}$ be the Lie superalgebra $\widetilde{p}(n)\,(n \geq 2)$ over $\mathbb{C}$, $U$ be its universal enveloping algebra. Recall that $\mathfrak{g}_0 \cong gl\,(n)$, and $\mathfrak{g}$ has a $\mathbb{Z}$-grading $\mathfrak{g} = \mathfrak{g}_{-1} \oplus \mathfrak{g}_0 \oplus \mathfrak{g}_1$, where $\mathfrak{g}_{-1}$ is isomorphic to $\Lambda^2 E^*$ and $\mathfrak{g}_1$ is isomorphic to $S^2 E$ as $\mathfrak{g}_0$-module, here $E$ stands for the standard module over $\mathfrak{g}_0$. Let $\mathfrak{g}_0' = [\mathfrak{g}_0, \mathfrak{g}_0] \cong sl\,(n)$. The grading of $\mathfrak{g}$ induces the grading of $U$ in the natural way. Thus

$$U = \bigoplus_{i=-n(n-1)/2}^{n(n+1)/2} U_i.$$

The standard matrix realization of $\widetilde{p}(n)$ is given by all block matrices

$$\begin{pmatrix} a & b \\ c & -a^t \end{pmatrix}$$

where $a$ is an arbitrary $n \times n$-matrix, $a^t$ denotes the transposed matrix, $b$ is a symmetric $n \times n$-matrix and $c$ is a skew-symmetric $n \times n$-matrix.

We fix a Cartan subalgebra $\mathfrak{h}$ in $\mathfrak{g}$, which coincides with Cartan subalgebra of $\mathfrak{g}_0$ and identify $\mathfrak{h}$ with $\mathfrak{h}^*$ by means of the $\mathfrak{g}_0$-invariant form $(A_1, A_2) = \operatorname{tr} A_1 A_2$ on $\mathfrak{g}_0$. If $\varepsilon_1, \ldots, \varepsilon_n$ is the standard orthogonal basis in $\mathfrak{h}$, then the roots of $\mathfrak{g}$ are of the following form:

Roots of $\mathfrak{g}_0$: $\varepsilon_i - \varepsilon_j$, $i \neq j$, $1 \leq i, j \leq n$

Roots of $\mathfrak{g}_{-1}$: $-\varepsilon_i - \varepsilon_j$, $1 \leq i < j \leq n$

Roots of $\mathfrak{g}_1$: $\varepsilon_i + \varepsilon_j$, $1 \leq i \leq j \leq n$

Each root has multiplicity one. For each root $\alpha$ we fix an element $X_\alpha$ from the root space such that $[X_\alpha, X_{-\alpha}] = H_\alpha$ with $(H_\alpha, H_\alpha) = 2$. If $-\alpha$ is not a root we choose $X_\alpha$ being an arbitrary non-zero element of the root space. For $\mu \in \mathfrak{h}^*$ we denote by $\mu_i$ its coordinate in the standard basis, i.e. $\mu = \sum_{i=1}^n \mu_i \varepsilon_i$.

We also denote by $(\cdot, \cdot)$ the bilinear symmetric form on $\mathfrak{h}^*$ induced by the form on $\mathfrak{h}$. One can check that $(\alpha, \beta) = \alpha\,(H_\beta)$ for any two roots $\alpha, \beta$ of $\mathfrak{g}_0$.

We choose a Borel subalgebra $\mathfrak{b} \supset \mathfrak{h}$ such that the positive roots are $\varepsilon_i - \varepsilon_j$ for $i < j$ and all roots of $\mathfrak{g}_1$. Let $\mathfrak{n} = [\mathfrak{b}, \mathfrak{b}]$. As was shown in [3], all irreducible finite dimensional $\mathfrak{g}$-modules are modules of highest weight. An irreducible module with highest weight $\mu$ is finite-dimensional if and only if $\mu_i - \mu_j$ is a non-negative integer for all $i < j$.

Finally let $\mathfrak{g}' = p(n) = [\mathfrak{g}, \mathfrak{g}]$, $\mathfrak{h}' = \mathfrak{h} \cap \mathfrak{g}'$. One can check that $\mathfrak{g}' = \mathfrak{g}_{-1} \oplus \mathfrak{g}'_0 \oplus \mathfrak{g}_1$. Let $U'$ denote the universal enveloping of $\mathfrak{g}'$.

## 3. Radical of $U$

In this section we will show that $U$ has a non-zero Jacobson radical and give some description of this radical. This radical was constructed in [4]. We put $d = n(n-1)/2$. Note that $d = \dim \mathfrak{g}_{-1}$.

Let $M_\mu^0$ denote Verma module over $\mathfrak{g}_0$ with highest weight $\mu$ and $M_\mu$ be Verma module over $\mathfrak{g}$. One can verify easily that

(3.1) $$M_\mu = \mathrm{Ind}_{\mathfrak{g}_0 \oplus \mathfrak{g}_1}^{\mathfrak{g}} M_\mu^0,$$

assuming $\mathfrak{g}_1 M_\mu^0 = 0$.

Let $v$ be a highest vector of $M_\mu$ and $X = \Pi_{i<j} X_{\varepsilon_i + \varepsilon_j}$, $Y = \Pi_{i<j} X_{-\varepsilon_i - \varepsilon_j}$. Finally let $\Delta(\mu) = \Pi_{i<j}(\mu_i - \mu_j + j - i - 1)$.

**Lemma 3.1.** $XYv = \Delta(\mu)v$.

*Proof.* We order the roots $\varepsilon_i + \varepsilon_j$ in the following way

$$\varepsilon_i + \varepsilon_j < \varepsilon_p + \varepsilon_q \text{ iff } i < p \text{ or } i = p \text{ and } j < q.$$

According to this order we index our roots $\alpha_1, \ldots, \alpha_d$.

Then whenever $\alpha < \beta < \gamma$ we have

$$[X_\alpha, X_{-\beta}] \in \mathfrak{b}_0,$$

$[[X_\alpha, X_{-\beta}], X_{-\gamma}]$ is either zero or $X_{-\delta}$, with $\gamma < \delta$.

Now let $w_k = X_{-\alpha_k} \ldots X_{-\alpha_d} v$, $w_{d+1} = v$. Using the above relations one verifies that

$$X_{\alpha_k} w_k = (\mu - \alpha_{k+1} - \cdots - \alpha_d, H_{\alpha_k}) w_{k+1}.$$

If $\alpha_k = \varepsilon_i + \varepsilon_j$, then $(\mu - \alpha_{k+1} - \cdots - \alpha_d, H_{\alpha_k}) = \mu_i - \mu_j + j - i - 1$. Repeating this argument one obtains $Xw_1 = \Delta(\mu)w_{d+1}$. □

Let $\mathfrak{h}_0^* = \{\mu \in \mathfrak{h}^* \mid \mu_i - \mu_j \notin \mathbb{Z}, i \neq j\}$. Clearly $\mathfrak{h}_0^*$ is a Zariski dense set in $\mathfrak{h}^*$.

**Lemma 3.2.** *If $\mu \in \mathfrak{h}_0^*$ then $M_\mu$ is irreducible.*

*Proof.* If we put $\deg v = 0$, then $M_\mu$ inherits $\mathbb{Z}$-grading from $U$, such that $M_\mu = (M_\mu)_0 + \cdots + (M_\mu)_{-d}$. Note that $(M_\mu)_{-d} = YM_\mu^0 = M_{\mu+\sigma}^0$, where $\sigma$ is the sum of all roots of $\mathfrak{g}_{-1}$. The highest vector of $M_{\mu+\sigma}^0$ is $Yv$. By (3.1) $M_\mu$ is free over $U(\mathfrak{g}_{-1})$. Therefore every non-zero submodule $N$ of $M_\mu$ has a nontrivial intersection with $M_{\mu+\sigma}^0$. The conditions on $\mu$ ensure that $M_{\mu+\sigma}^0$ is irreducible. Thus, $N$ must contain $Yv$, and therefore $XYv = \Delta(\mu)v$. As $\mu \in \mathfrak{h}_0^*$, $\Delta(\mu) \neq 0$, we have $v \in N$ and hence $N = M_\mu$. □

**Lemma 3.3.** *Let* $P_\mu = \operatorname{Ind}_{\mathfrak{g}_0}^{\mathfrak{g}} M_\mu^0$.

*(a) If* $\mu \in \mathfrak{h}_0^*$, *then* $P_\mu$ *has finite length with irreducible subquotients* $M_{\mu+\gamma}$, *where* $\gamma$ *runs over sums of non-repeating roots of* $\mathfrak{g}_1$. *Note that the length* $s$ *of* $P_\mu$ *is equal to* $2^{n(n+1)/2}$.

*(b)* $\bigoplus_{\mu \in \mathfrak{h}_0^*} P_\mu$ *is a faithful* $U$-*module.*

*Proof.* To prove (a) use

$$P_\mu = \operatorname{Ind}_{\mathfrak{g}_0 \oplus \mathfrak{g}_1}^{\mathfrak{g}} \left( \operatorname{Ind}_{\mathfrak{g}_0}^{\mathfrak{g}_0 \oplus \mathfrak{g}_1} M_\mu^0 \right).$$

A filtration of $\operatorname{Ind}_{\mathfrak{g}_0}^{\mathfrak{g}_0 \oplus \mathfrak{g}_1} M_\mu^0$ with $\mathfrak{g}_0$-irreducible quotients $M_{\mu+\gamma}^0$ induces the filtration on $P_\mu$ with irreducible quotients $M_{\mu+\gamma}$.

To show (b) recall that $\bigoplus_{\mu \in \mathfrak{h}_0^*} M_\mu^0$ is a faithful $U(\mathfrak{g}_0)$-module. Therefore $\bigoplus_{\mu \in \mathfrak{h}_0^*} P_\mu = \operatorname{Ind}_{\mathfrak{g}_0}^{\mathfrak{g}} \left( \oplus_{\mu \in \mathfrak{h}_0^*} M_\mu^0 \right)$ is a faithful $U$-module. $\square$

The following theorem is proved in [4]. If $M$ is a module over an associative algebra $\mathcal{A}$, then $\operatorname{Ann}_{\mathcal{A}} M$ denotes the annihilator of $M$ in $\mathcal{A}$.

**Theorem 3.4.** *The ideal* $I = \cap_{\mu \in \mathfrak{h}_{\mu} \in \mathfrak{h}} \operatorname{Ann}_U M_\mu = \cap_{\mu \in \mathfrak{h}_0^*} \operatorname{Ann}_U M_\mu$ *coincides with the radical of* $U$.

*Proof.* The assertion $\cap_{\mu \in \mathfrak{h}^*} \operatorname{Ann}_U M_\mu = \cap_{\mu \in \mathfrak{h}_0^*} \operatorname{Ann}_U M_\mu$ follows from Zariski density of $\mathfrak{h}_0^*$. By Lemma 3.2 $M_\mu$ is irreducible if $\mu \in \mathfrak{h}_0^*$, and therefore $\operatorname{rad} U \subset I$. To show that $I \subset \operatorname{rad} U$ it suffices to prove that $I$ is nilpotent. Note that $I^s(P_\mu) = 0$ for all $\mu \in \mathfrak{h}_0^*$. Since $\oplus_{\mu \in \mathfrak{h}_0^*} P_\mu$ is faithful we obtain $I^s = 0$. $\square$

**Corollary 3.5.** *Let* $\bar{U} = U/I$. *Then the* $\mathbb{Z}$-*grading on* $U$ *induces the* $\mathbb{Z}$-*grading on* $\bar{U}$ *such that* $\bar{U} = \bar{U}_{-d} \oplus \cdots \oplus \bar{U}_d$. *The natural maps* $U(\mathfrak{g}_0 \oplus \mathfrak{g}_{-1}) \to \bar{U}$ *and* $\mathfrak{g} \to \bar{U}$ *are injective. Finally* $\cap_{\mu \in \mathfrak{h}_0^*} \operatorname{Ann}_{\bar{U}} M_\mu = \{0\}$.

The following theorem can be proven in the same manner as Theorem 3.4.

**Theorem 3.6.** *The ideal* $I' = \cap_{\mu \in (\mathfrak{h}')^*} \operatorname{Ann}_{U'} M_\mu$ *coincides with the radical of* $U'$.

## 4. THE CENTER OF $\bar{U}$

In this section we will describe the center $Z$ of $\bar{U}$. As follows from Corollary 3.5 $U(\mathfrak{h})$ is a subalgebra of $\bar{U}$. Theorem 3.4 implies that $I$ annihilates highest vectors of all Verma modules. Therefore $I \subset U\mathfrak{n}^+$. Hence the standard Harish-Chandra projection $h : U \to U(\mathfrak{h})$ with the kernel $\mathfrak{n}^- U + U\mathfrak{n}^+$ induces the map $\bar{U} \to S(\mathfrak{h})$. Thus, one can define the Harish-Chandra homomorphism $h : Z \to U(\mathfrak{h}) = S(\mathfrak{h})$. We identify $S(\mathfrak{h})$ with polynomial algebra on $\mathfrak{h}^*$. As in the usual case for any $z \in Z$

$$(4.1) \qquad\qquad z|_{M_\mu} = h(z)(\mu) \operatorname{Id}.$$

Therefore the last statement of Corollary 3.5 implies that $h$ is injective. We will describe $h(Z)$.

Let $W$ be the Weyl group of $\mathfrak{g}_0$, note that $W$ is isomorphic to $S_n$. Define a $W$-action on $\mathfrak{h}^*$ by $\mu^w = w(\mu + \rho_0) - \rho_0$, where $\rho_0$ is the half sum of even positive roots. Denote by $S(\mathfrak{h})^W$ the subring of $W$-invariant polynomials on $\mathfrak{h}^*$. Let $\Theta(\mu) = \Pi_{i \neq j}((\mu + \rho_0, \varepsilon_i - \varepsilon_j) - 1)$. Note that $\Theta \in S(\mathfrak{h})^W$.

**Lemma 4.1.** $h(Z) \subset S(\mathfrak{h})^W$.

*Proof.* Let $\mathfrak{d} = \{\mu \in \mathfrak{h}^* \mid \mu_i - \mu_j \in \mathbb{Z}_{\geq 0}, i < j\}$. Then $\mathfrak{d}$ is a Zariski dense set in $\mathfrak{h}^*$. Furthermore for any $w \in W$

$$\mathrm{Hom}_{\mathfrak{g}}(M_{\mu^w}, M_\mu) = \mathrm{Hom}_{\mathfrak{g}_0}(M_{\mu^w}^0, M_\mu^0) = \mathbb{C}.$$

Therefore every element $z \in Z$ acts by the same scalar on $M_{\mu^w}$ and on $M_\mu$. By (4.1) every $p \in h(Z)$ satisfies $p(\mu^w) = p(\mu)$ for all $\mu \in \mathfrak{d}$. Since $\mathfrak{d}$ is Zariski dense, this implies that $p(\mu^w) = p(\mu)$ for all $\mu \in \mathfrak{h}^*$. $\square$

**Lemma 4.2.** *Let* $\mu \in \mathfrak{h}^*$ *be such that* $\mu_{n-1} = \mu_n$ *and* $p \in h(Z)$. *Then* $p(\mu + t(\varepsilon_{n-1} + \varepsilon_n)) = p(\mu)$ *for any* $t \in \mathbb{C}$.

*Proof.* Let $v$ be a highest vector of Verma module $M_\mu$, $u = X_{\varepsilon_n - \varepsilon_{n-1}}v$ and $w = X_{-\varepsilon_n - \varepsilon_{n-1}}v$. An easy calculation shows that $u$ is $\mathfrak{n}^+$-invariant, and $\mathfrak{n}^+ w = \mathbb{C}u$. Therefore $N = \bar{U}u$ is a proper submodule in $M_\mu$ and $\mathrm{Hom}_{\mathfrak{g}}(M_{\mu - \varepsilon_{n-1} - \varepsilon_n}, M_\mu/N) \neq 0$. By (4.1) this implies $p(\mu) = p(\mu - \varepsilon_{n-1} - \varepsilon_n)$. Applying this relation several times we obtain $p(\mu) = p(\mu + t(\varepsilon_{n-1} + \varepsilon_n))$ for all $t \in \mathbb{Z}$, which by Zariski density of $\mathbb{Z}$ in $\mathbb{C}$ implies the same relation for all $t \in \mathbb{C}$. $\square$

**Corollary 4.3.** *Let* $(\mu + \rho_0, \varepsilon_i - \varepsilon_j) = 1$ *and* $p \in h(Z)$. *Then* $p(\mu) = p(\mu + t(\varepsilon_i + \varepsilon_j))$ *for all* $t \in \mathbb{C}$.

**Lemma 4.4.** $h(Z) \subset \Theta S(\mathfrak{h})^W + \mathbb{C}$.

*Proof.* Let

$$l_{i,j} = \{\mu \in \mathfrak{h}^* \mid (\mu + \rho_0, \varepsilon_i - \varepsilon_j) = 1\}.$$

We want to show that if $p \in h(Z)$, then $p$ is constant on $\cup_{i \neq j} l_{i,j}$. Since $p$ is $W$-invariant, it suffices to show that $p$ is constant on one hyperplane $l_{i,j}$, for example on $l_{1,2}$. If $n = 2$, the statement is trivial. Let $\mu \in l_{1,2}$, i.e. $\mu_1 = \mu_2$. Put $\lambda = \mu + (\mu_3 - \mu_1)(\varepsilon_1 + \varepsilon_2)$. By Corollary 4.3 $p(\mu) = p(\lambda)$. Note that $\lambda \in l_{2,3}$. Put $\nu = \lambda + (\mu_1 - \mu_3 + 1)(\varepsilon_2 + \varepsilon_3)$. Then $p(\mu) = p(\lambda) = p(\nu)$. Let $w = (132) \in W$. Then $\nu^w = \mu + 2\varepsilon_3$. Thus we obtain the relation $p(\mu) = p(\mu + 2\varepsilon_3)$. In the same way as in the proof of lemma 5 we can get that $p(\mu) = p(\mu + t\varepsilon_3)$ for all $t \in \mathbb{C}$. Using $W$-invariance of $p$ we obtain that $p(\mu) = p(\mu + t\varepsilon_i)$ for all $t \in \mathbb{C}$ and $i = 3, \ldots, n$. This implies that $p$ is constant on $l_{1,2}$ and therefore $p$ is constant on $\cup_{i \neq j} l_{i,j}$. $\square$

For any $x \in U$ denote by $x_0$ the image of $x$ under the projection onto $U(\mathfrak{g}_0)$ with the kernel $\mathfrak{g}_1 U + U\mathfrak{g}_{-1}$. Recall that $U = U(\mathfrak{g}_1) \otimes U(\mathfrak{g}_0) \otimes U(\mathfrak{g}_{-1})$ and $U(\mathfrak{g}_1) = S(\mathfrak{g}_1) \subset U$.

**Lemma 4.5.** $(YS^d(\mathfrak{g}_1))_0 \neq 0$.

*Proof.* We will show that $(YX)_0 \neq 0$. Consider the Verma module $M_\mu$ with a highest vector $v$. Then $(YX)_0 Yv = YXYv = \Delta(\mu) Yv$ by Lemma 3.1. Therefore $(YX)_0 \neq 0$. □

Let $U^k$ be the $k$-th term of the natural filtration of $U$. Lemma 4.5 enables us to define a non-trivial $\mathrm{ad}_{\mathfrak{g}_0'}$-invariant map

$$\varphi : S^d(\mathfrak{g}_1) \to U^d(\mathfrak{g}_0)$$

by

$$\varphi(x) = (Yx)_0.$$

Let $\{X_i\}$ be a basis in some $\mathrm{ad}_{\mathfrak{g}_0}$-invariant subspace of $S^d(\mathfrak{g}_1)$ complementary to $\mathrm{Ker}\,\varphi$. Then $\{(YX_i)_0\}$ is a linearly independent set in $U^d(\mathfrak{g}_0)$ and could be extended to some basis of $U^d(\mathfrak{g}_0)$. Recall that $U^d(\mathfrak{g}_0) \cong S^{\leq d}(\mathfrak{g}_0)$ as $\mathrm{ad}_{\mathfrak{g}_0}$ module and it admits a non-degenerate $\mathrm{ad}_{\mathfrak{g}_0}$-invariant bilinear symmetric form $B$. Use the notation $(YX_i)_0'$ for the vectors of the dual basis with respect to $B$. Let

$$T = \Sigma X_i (YX_i)_0'.$$

Then $T$ is $\mathrm{ad}_{\mathfrak{g}_0'}$ invariant and

$$(4.2) \qquad (YT)_0 = \Sigma \left(YX_i (YX_i)_0'\right)_0 = \Sigma (YX_i)_0 (YX_i)_0' \neq 0.$$

**Lemma 4.6.** *(a) Let $Y^{\mathrm{ad}} = \Pi_{i<j}\,\mathrm{ad}_{X_{-\varepsilon_i-\varepsilon_j}}$ and $S = Y^{\mathrm{ad}}(T)$, and $\bar{S}$ be its image in $\bar{U}$. Then $\bar{S} \in Z$ and $h\left(\bar{S}\right)$ is proportional to $\Theta$.*

*(b) Let $z \in Z(U(\mathfrak{g}_0))$. Put $S_z = Y^{\mathrm{ad}}(zT)$, and let $\bar{S}_z$ be its image in $\bar{U}$. Then $\bar{S}_z \in Z$ and $h\left(\bar{S}_z\right) = h\left(\bar{S}\right) h(z)$.*

*Proof.* (a) First of all $\bar{S} \neq 0$. Indeed, if $v$ is a highest vector of $M(\mu)$ then

$$\bar{S}Yv = SYv = Y^{\mathrm{ad}}(T) Yv = YTYv = (YT)_0 Yv.$$

Since $(YT)_0 \neq 0$, $(YT)_0 Yv \neq 0$ for some $\mu$, and therefore $\bar{S} \neq 0$.

Note that $\bar{T} \in \bar{U}_d$, hence $\mathfrak{g}_1 \bar{T} = \bar{T}\mathfrak{g}_1 = 0$, and $\mathrm{ad}_{\mathfrak{g}_0'} \bar{T} = 0$. Thus, $\bar{T}$ is $\mathrm{ad}_{\mathfrak{b}}$-invariant vector of weight $-\sigma$. Let $V$ be $\mathrm{ad}_{\mathfrak{g}}$-submodule of $\bar{U}$ generated by $\bar{T}$. Since $\bar{S} \neq 0$, $V \cong \mathrm{Ind}_{\mathfrak{g}_0 \oplus \mathfrak{g}_1}^{\mathfrak{g}} \mathbb{C}\bar{T}$ and $\bar{S}$ generates the trivial $\mathrm{ad}_{\mathfrak{g}}$-submodule of $V$. Thus, $\bar{S} \in Z$.

To find $h\left(\bar{S}\right)$ note that $\bar{S}$ acts by zero on the trivial $\mathfrak{g}$-module, hence $h\left(\bar{S}\right)(0) = 0$. On the other hand the degree of the polynomial $h\left(\bar{S}\right)$ is not greater than $2d$ as $\bar{S} \in \bar{U}^{2d}$. Therefore, by Lemma 4.4, $h\left(\bar{S}\right)$ must be proportional to $\Theta$.

(b) Arguments similar to (a) show that $\bar{S}_z \in Z$. Finally, if $v$ is a highest vector of $M_\mu$, then

$$(4.3) \qquad h\left(\bar{S}_z\right)(\mu) v = \bar{S}_z v = Y^{\mathrm{ad}}(zT) v = zTYv = zY^{\mathrm{ad}}(T) v = zSv.$$

Note that $zSv = (h(z))(\mu) Sv$. Therefore $h\left(\bar{S}_z\right) = h\left(\bar{S}\right) h(z)$. □

**Corollary 4.7.** $\mathbb{C} + \Theta S(\mathfrak{h})^W \subset h(Z)$.

Lemma 4.4 and Corollary 4.7 immediately imply

**Theorem 4.8.** $h(Z) = \mathbb{C} + \Theta S(\mathfrak{h})^W$.

Let $\bar{U}' = U'/I'$, $Z'$ be the center of $\bar{U}'$ and $h \colon Z' \to S(\mathfrak{h}')^*$ be the Harish-Chandra homomorphism. One can prove the following theorem by repeating the arguments of this section.

**Theorem 4.9.** $h(Z') = \mathbb{C} + \Theta S(\mathfrak{h}')^W$.

## 5. The category of $\bar{U}$-modules

Let $\mu \in \mathfrak{h}^*$. Define $\chi_\mu \colon Z \to \mathbb{C}$ by the formula $\chi_\mu(z) = h(z)(\mu)$ for all $z \in Z$. We call $\mu$ *typical* if $\Theta(\mu) \neq 0$. Otherwise we call $\mu$ *atypical*. The following statement follows immediately from Theorem 4.8.

**Lemma 5.1.** *Let* $\mu, \nu \in \mathfrak{h}^*$. *Then* $\chi_\mu = \chi_\nu$ *iff either both* $\mu$ *and* $\nu$ *are atypical or* $\mu$ *and* $\nu$ *are both typical and* $\mu = \nu^w$ *for some* $w \in W$

Define $\bar{U}_\mu = \bar{U}/\bar{U} \operatorname{Ker} \chi_\mu$. Let $U^0_\mu = U(\mathfrak{g}_0)/U(\mathfrak{g}_0) \operatorname{Ker} \chi^0_\mu$, where $\chi^0_\mu \colon Z(U(\mathfrak{g}_0)) \to \mathbb{C}$ is defined using the Harish-Chandra projection restriction to $Z(U(\mathfrak{g}_0))$. By $\bar{U}_\mu -$ mod (respectively $U^0_\mu -$ mod ) we denote the category of all left $\bar{U}_\mu$-modules (resp. $U^0_\mu$-modules). In the similar way one defines $\bar{U}'_\mu$ for $\mu \in (\mathfrak{h}')^*$.

We will prove the following

**Theorem 5.2.** *If* $\mu \in \mathfrak{h}^*$ *is typical then the categories* $\bar{U}_\mu -$ mod *and* $U^0_\mu -$ mod *are equivalent.*

Let us introduce the following notations. Put $z_0 = \operatorname{diag}(1_n, -1_n) \in Z(\mathfrak{g}_0)$. We assume now that $\mu$ is typical and fixed. If $M$ is a $\mathfrak{g}$-module we put

$$M^{\mathfrak{g}_1} = \{m \in M \mid \mathfrak{g}_1 m = 0\}$$
$$M_0 = \{m \in M \mid z_0 m = \mu(z_0) m\}.$$

**Lemma 5.3.** *If* $M$ *is an* $\bar{U}_\mu$-module then
   (a) $M^{\mathfrak{g}_1}$ *is an* $U^0_\mu$-module;
   (b) $M^{\mathfrak{g}_1} = M_0$;
   (c) $z_0$ *acts semi-simply on any* $\bar{U}_\mu$-module;
   (d) The functor $\operatorname{Inv} \colon \bar{U}_\mu -$ mod $\to U^0_\mu -$ mod , defined by $\operatorname{Inv} M = M_0 = M^{\mathfrak{g}_1}$ is faithful and exact.*

*Proof.* The first statement follows from the formula $\bar{S}_z m = z \bar{S} m$ for any $m \in M^{\mathfrak{g}_1}$, $z \in Z(U(\mathfrak{g}_0))$. Since $\bar{S} m = \chi_\mu(\bar{S}) m \neq 0$, we get

$$zm = \frac{\bar{S}_z m}{\chi_\mu(\bar{S})} = \chi^0_\mu(z) m.$$

By (a) $z_0 m = \mu(z_0) m$ for any $m \in M^{\mathfrak{g}_1}$, and therefore $M^{\mathfrak{g}_1} \subseteq M_0$. Let $m \in M_0$. Then $M' = U(\mathfrak{g}_1) m$ is a $z_0$-invariant finite-dimensional subspace of $M$ with $z_0$ weight decomposition $M' = \oplus_{i \geq 0} M'_{\mu(z_0)+i}$. In particular $M'_{\mu(z_0)} = \mathbb{C}m$. Since $(M')^{\mathfrak{g}_1} \neq 0$, and $(M')^{\mathfrak{g}_1} \subseteq M_0$, we obtain $m \in (M')^{\mathfrak{g}_1}$ and $M' = \mathbb{C}m$. Therefore $M^{\mathfrak{g}_1} = M_0$.

Faithfulness of Inv follows from the fact that $U(\mathfrak{g}_1)$ is isomorphic to finite-dimensional Grassmann algebra.

To prove $(c)$ let $N$ be a maximal submodule of $M$ on which $z_0$ acts semi-simply. By $(b)$ $N^{\mathfrak{g}_1} = M^{\mathfrak{g}_1}$. Assume that $N \neq M$. Then $(M/N)^{\mathfrak{g}_1} \neq 0$. Let $M' \subset M$ be the preimage of $(M/N)^{\mathfrak{g}_1}$ under the natural projection $M \to M/N$. One can find $v \in M' \backslash N$ such that $z_0 v = \mu(z_0)v + w$ and $z_0 w = \mu(z_0)w$. By $(b)$ $w \in M^{\mathfrak{g}_1}$. Therefore $yw = 0$ and $yz_0 v = \mu(z_0)yv$ for any $y \in \mathfrak{g}_1$. But $z_0 yv = yz_0 v + yv$, and therefore $z_0 yv = (\mu(z_0) + 1)yv$. On the other hand, there exist $y_1, \ldots, y_k \in \mathfrak{g}_1$ such that $y_1 \ldots y_k yv \in M^{\mathfrak{g}_1} = M_0$. Therefore $z_0 yv = (\mu(z_0) - k)yv$. Hence $yv = 0$ for any $y \in \mathfrak{g}_1$, but that implies $z_0 v = \mu(z_0)v$. Therefore $v \in N$. Contradiction.

Exactness is a direct corollary of the identity Inv $M = M_0$. $\qquad \square$

Let $\bar{U}^+$ denote the image of $U(\mathfrak{g}_0 \oplus \mathfrak{g}_1)$ in $\bar{U}$ under the natural projection. Since $U(\mathfrak{g}_0)$ is a subalgebra of $\bar{U}$ it is possible to define a $\bar{U}^+$-module structure on every $\mathfrak{g}_0$-module by putting $\mathfrak{g}_1 N = 0$.

Let Ind $N = \bar{U} \otimes_{\bar{U}^+} N$.

**Lemma 5.4.** *Let $N$ be a $U_\mu^0$-module. Then*

(a) Ind $N$ *is $\bar{U}_\mu$-module.*

(b) Ind $N \cong U(\mathfrak{g}_{-1}) \otimes_{\mathbb{C}} N$ *as $\mathfrak{g}_0$-module.*

(c) *The functor* Ind$: U_\mu^0 - \text{mod} \to \bar{U}_\mu - \text{mod}$ *is exact and faithful.*

(d) Inv (Ind $N$) $\cong N$.

*Proof.* Since Ind $N = \bar{U}N$, it suffices to check that $zv = \chi_\mu(z)v$ for every $z \in Z$ and $v \in N$. Define the projection $h_0 : \bar{U}_0 \to U(\mathfrak{g}_0)$ with the kernel $\bar{U}_0 \cap \bar{U}\mathfrak{g}_1$. One can easily verify that $h_0(Z) \subset Z(U(\mathfrak{g}_0))$, $zv = h_0(z)v$ and $\chi_\mu(z) = \chi_\mu^0(h_0(z))$. Now (a) follows immediately.

To prove (b) recall that $U(\mathfrak{g}_{-1})$ is a subalgebra of $\bar{U}$ and Ind $N$ is free over $U(\mathfrak{g}_{-1})$.

(c) follows from (b).

To prove (d) use that Ind $N = N \oplus \mathfrak{g}_{-1}U(\mathfrak{g}_{-1})$ Ind $N$. Let

$$N' := \text{Inv}(\text{Ind } N) \cap (\mathfrak{g}_{-1} \text{Ind } N).$$

It suffices to show that $N' = \{0\}$. Indeed, let $M' \subset $ Ind $N$ be a submodule generated by $N'$. Then $M' \subseteq \mathfrak{g}_{-1}$ Ind $N$, and therefore $TM' = 0$. But then $\bar{S}M' = 0$. Since $\mu$ is typical, $\bar{S}$ acts on $M'$ by a non-zero constant. Hence $M' = \{0\}$ and $N' = \{0\}$. $\quad \square$

**Lemma 5.5.** *The functors* Inv$: U_\mu^0 - \text{mod} \to \bar{U}_\mu - \text{mod}$ *and* Ind$: \bar{U}_\mu - \text{mod} \to U_\mu^0 - \text{mod}$ *establish equivalence of categories.*

*Proof.* By Lemma 5.4 $(d)$ Inv (Ind $N$) $\cong N$. Therefore it suffices to show that Ind (Inv $M$) $\cong M$. Since Inv (Ind (Inv $M$)) $\cong$ Inv $M$ and Inv is faithful and exact, we obtain Ind (Inv $M$) $\cong M$. $\qquad \square$

Note that Lemma 5.5 implies Theorem 5.2.

*Remark* 5.6. Note that in the same way one can prove that the category of right $\bar{U}_\mu$-modules and the category of right $U^0_{\mu+\sigma}$-modules are equivalent. The shift by $\sigma$ comes from rewriting formula (4.3) for right $\mathfrak{g}$-modules.

**Theorem 5.7.** *For a typical $\mu \in \mathfrak{h}^*$ the algebra $\bar{U}_\mu$ is isomorphic to the matrix algebra $U^0_\mu \otimes_{\mathbb{C}} \operatorname{End}_{\mathbb{C}}(U(\mathfrak{g}_{-1}))$.*

*Proof.* We are using the following classical result (see [2]). Let $A$ and $B$ be $\mathfrak{g}_0$-modules. Denote by $F(A, B)$ the subspace of $\operatorname{Hom}_{\mathbb{C}}(A, B)$ on which the adjoint action of $\mathfrak{g}_0$ is locally finite. The natural homomorphism $\rho_0 : U^0_\mu \to F(M^0_\mu, M^0_\mu)$ is an isomorphism. We will show that the natural homomorphism $\rho : \bar{U}_\mu \to F(M_\mu, M_\mu)$ is also an isomorphism. First of all since $M_\mu = \operatorname{Ind} M^0_\mu$ and $M^0_\mu$ is faithful, Morita equivalence of $U^0_\mu$ and $\bar{U}_\mu$ implies that $M_\mu$ must be faithful. Thus $\rho$ is injective. To show the surjectivity of $\rho$ consider $F(M_\mu, M_\mu)$ as a left $\bar{U}_\mu$-module. Then by Theorem 5.2

$$F(M_\mu, M_\mu) = \operatorname{Ind}(\operatorname{Inv} F(M_\mu, M_\mu)) = \operatorname{Ind} F(M_\mu, M^0_\mu).$$

Note that $F(M_\mu, M^0_\mu)$ is a right $\bar{U}_\mu$-module and therefore by Remark 5.6

$$F(M_\mu, M^0_\mu) = \operatorname{Ind}(\operatorname{Inv} F(M_\mu, M^0_\mu)) = \operatorname{Ind} F(M^0_{\mu+\sigma}, M^0_\mu).$$

Thus, it suffices to check that $F(M^0_{\mu+\sigma}, M^0_\mu) \subseteq \operatorname{Im} \rho$. Indeed, it follows from

$$F(M^0_{\mu+\sigma}, M^0_\mu) = F(M^0_\mu, M^0_\mu) \rho(T) = \rho(U(\mathfrak{g}_0) T).$$

Finally, the isomorphism of $\mathfrak{g}_0$-modules $M_\mu \cong M^0_\mu \otimes_{\mathbb{C}} U(\mathfrak{g}_{-1})$ implies

$$F(M_\mu, M_\mu) \cong F(M^0_\mu, M^0_\mu) \otimes_{\mathbb{C}} \operatorname{End}_{\mathbb{C}}(U(\mathfrak{g}_{-1})),$$

and that completes the proof of the theorem. $\qquad\square$

**Corollary 5.8.** *Let $M_0$ be an irreducible $U^0_\mu$-module for some typical $\mu \in \mathfrak{h}^*$. Then $M = \operatorname{Ind}^{\mathfrak{g}}_{\mathfrak{g}_0 \oplus \mathfrak{g}_1} M_0$ is an irreducible $\mathfrak{g}$-module. Moreover, if $M_0$ is a semi-simple $\mathfrak{h}$-module with finite weight multiplicities, then $M$ is also $\mathfrak{h}$-semi-simple with finite weight multiplicities and $\operatorname{ch} M = \operatorname{ch} M_0 \Pi_{i \neq j}(1 + \varepsilon e^{-\varepsilon_i - \varepsilon_j})$.*

**Theorem 5.9.** *Let $\nu \in (\mathfrak{h}')^*$ and $\mu \in \mathfrak{h}^*$ be such that $\mu_{|\mathfrak{h}'} = \nu$. Then the algebra $\bar{U}'_\nu$ is isomorphic to $\bar{U}_\mu$.*

*Proof.* Let $(U^0_\nu)' = U(\mathfrak{g}'_0)/U(\mathfrak{g}'_0) \operatorname{Ker} \chi^0_\nu$. One can prove repeating the above arguments that $\bar{U}'_\nu$ is isomorphic to $(U^0_\nu)' \otimes \operatorname{End}_{\mathbb{C}}(U(\mathfrak{g}_{-1}))$. Now the theorem follows from the obvious isomorphism $(U^0_\nu)' \cong U^0_\mu$. $\qquad\square$

## 6. Geometric realization of $\operatorname{Gr} \bar{U}$

As in the classical case, $\operatorname{Gr} U \cong S(\mathfrak{g})$ is a supercommutative algebra with a Poisson bracket. It can be considered as the algebra of polynomials on $\mathfrak{g}^*$. As $I$ is a two-sided ideal in $U$, $\operatorname{Gr} I$ is a Poisson ideal, and the quotient algebra $\operatorname{Gr} \bar{U} = \operatorname{Gr} U / \operatorname{Gr} I$ is also a supercommutative algebra with a Poisson bracket. Being a quotient algebra of $\operatorname{Gr} U$, it can be considered as an algebra of functions on an appropriate closed Poisson subvariety (or a subscheme) of $\mathfrak{g}^*$. The reader can find all necessary details on supergeometry in [5] and [6]. Below we repeat only the principal notions.

Let $\mathcal{V}$ be a vector superspace. Then it is a supermanifold with the supercommutative ring of regular function $\mathcal{O}(\mathcal{V}) = S(\mathcal{V}^*)$. A $\mathbb{Z}_2$-graded ideal $I_{\mathcal{M}} \subset \mathcal{O}(\mathcal{V})$ defines an affine superscheme $\mathcal{M}$ inside $\mathcal{V}$. The supercommutative ring $\mathcal{O}(\mathcal{M}) = \mathcal{O}(\mathcal{V})/I_{\mathcal{M}}$ is called the ring of regular functions on $\mathcal{M}$. The ideal $J_{\mathcal{M}} \subset \mathcal{O}(\mathcal{M})$ generated by all odd elements of $\mathcal{O}(\mathcal{M})$ defines the affine scheme $\mathcal{M}_0 \subset \mathcal{V}_0$ which is called the underlying scheme of $\mathcal{M}$. We say that $\mathcal{M}$ is an affine supervariety if $\mathcal{M}_0$ is an affine variety, i.e. $\mathcal{O}(\mathcal{M}_0)$ does not have nilpotents.

Consider $\mathfrak{g}^*$ as a supermanifold with the coadjoint action of the supergroup $G$. Its underlying manifold is $\mathfrak{g}_0^*$. Let $\xi \colon G \times \mathfrak{g}_0^* \to \mathfrak{g}^*$ be the morphism of supermanifolds induced by the action of $G$ on $\mathfrak{g}^*$ and the canonical embedding $\mathfrak{g}_0^* \subset \mathfrak{g}^*$. Let $I_Q = \operatorname{Ker} \xi^*$ and $Q$ be the closed subscheme of $\mathfrak{g}^*$ defined by the ideal $I_Q$. Since the underlying variety of $Q$ coincides with $\mathfrak{g}_0^*$, $Q$ is an affine supervariety. By analogy with even situation one can say that $Q$ is the closure of $G\mathfrak{g}_0^*$. Note that $I_Q$ is $\operatorname{ad}_{\mathfrak{g}}$-invariant, hence $\{x, I_Q\} \subset I_Q$ for all $x \in \mathfrak{g}$. Therefore $I_Q$ is a Poisson ideal in $\mathcal{O}(\mathfrak{g}^*)$ and $Q$ is a (singular) Poisson subvariety of $\mathfrak{g}^*$.

Let us describe the geometry of $Q$. Denote by $R_0$ the open subset of all regular elements of $\mathfrak{g}_0^*$.

**Lemma 6.1.** *For $x \in \mathfrak{g}_0^*$, let $\mathfrak{g}_x$ be the Lie superalgebra of the stabilizer of $x$. Then $\dim(\mathfrak{g}_x \cap \mathfrak{g}_1) \geq n$. If $x \in R_0$ then $\dim(\mathfrak{g}_x \cap \mathfrak{g}_1) = n$ and $\mathfrak{g}_x \cap \mathfrak{g}_{-1} = \{0\}$. If $x \notin R_0$, then $\dim(\mathfrak{g}_x \cap \mathfrak{g}_1) > n$ and $\dim(\mathfrak{g}_x \cap \mathfrak{g}_{-1}) > 0$.*

*Proof.* Although in our case the adjoint and coadjoint representations are not isomorphic, there is a convenient matrix realization of $\mathfrak{g}^*$. We identify $\mathfrak{g}^*$ with the space of all matrices of the form

$$(6.1) \qquad \begin{pmatrix} q & s \\ u & q^t \end{pmatrix}$$

where $q$ is an arbitrary $n \times n$-matrix, $s$ is a skew-symmetric $n \times n$ matrix and $u$ is a symmetric $n \times n$-matrix. The pairing between $\mathfrak{g}$ and $\mathfrak{g}^*$ is given by the form $\langle x, y \rangle =$str $xy$, and the action of $\mathfrak{g}$ on $\mathfrak{g}^*$ is given by the supercommutator.

Let $x = \begin{pmatrix} q & 0 \\ 0 & q^t \end{pmatrix}$. Then

$$\mathfrak{g}_x \cap \mathfrak{g}_{-1} = \left\{ c \in \mathfrak{g}_{-1} \mid cq - q^t c = 0 \right\}$$

and in the same way

$$\mathfrak{g}_x \cap \mathfrak{g}_1 = \left\{ b \in \mathfrak{g}_1 \mid qb - bq^t = 0 \right\}.$$

Since $b$ is symmetric and $qb - bq^t$ is skewsymmetric,

$$\dim \mathfrak{g}_x \cap \mathfrak{g}_1 \geq \dim \mathfrak{g}_1 - \dim \mathfrak{g}_{-1} = n.$$

That proves the first statement.

To check two other statements note that $\dim \mathfrak{g}_x \cap \mathfrak{g}_{-1}$ and $\dim \mathfrak{g}_x \cap \mathfrak{g}_1$ are constant on the conjugacy class of $q$. Therefore it is sufficient to prove the statements for a canonical Jordan form of $q$ which is a straightforward calculation. $\square$

**Corollary 6.2.** $\dim Q = (n^2, n^2 - n)$.

The next step is to construct "a resolution of singularities" of $Q$.

Let $\mathfrak{p} = \mathfrak{g}_0 + \mathfrak{g}_1$, $P$ be the subgroup with the Lie superalgebra $\mathfrak{p}$, and $K = G/P$. Note that $P$ is the normalizer of $\mathfrak{g}_1$ and therefore $K$ is the supermanifold of all Lie subalgebras of $\mathfrak{g}$ $\mathrm{Ad}_G$-conjugate to $\mathfrak{g}_1$. Let $G_{-1}$ be the supergroup with the Lie superalgebra $\mathfrak{g}_{-1}$. Since $G_{-1}$ is supercommutative, the exponential map defines an isomorphism of algebraic supermanifolds $G_{-1}$ and $\mathfrak{g}_{-1}$. Since $G_{-1}$ acts simply transitively on $K$, $K$ is isomorphic to $\mathfrak{g}_{-1}$ as an algebraic supermanifold, and $\mathcal{O}(K) \cong S\left(\mathfrak{g}_{-1}^*\right)$.

Let $V$ be the vector bundle on $K$ induced by the $P$-module $(\mathfrak{g}/\mathfrak{g}_1)^*$. A point of $V$ is a pair $(\mathfrak{g}_1', \phi)$, where $\mathfrak{g}_1'$ is a subalgebra $\mathrm{Ad}_G$-conjugate to $\mathfrak{g}_1$ and $\phi$ is a linear functional on $\mathfrak{g}$ annihilating $\mathfrak{g}_1'$. (Here points are in the sense of the "functor of points" see [5]).

Define the map $r : V \to \mathfrak{g}^*$ by putting $r\left(\mathfrak{g}_1', \phi\right) = \phi$.

Note that $\mathfrak{g}_{-1}^* \times \mathfrak{g}_0^*$ is embedded in $V$ as the fiber over the point (in usual sense) $[\mathfrak{g}_1] \in G/P$. Since the base $K$ of the bundle $V$ is identified with $\mathfrak{g}_{-1}$, one can identify $V$ with $\mathfrak{g}_{-1} \times \mathfrak{g}_{-1}^* \times \mathfrak{g}_0^*$ using the action of $G_{-1}$. So a point $v \in V$ is a triple $v = (c, d, y) \in \mathfrak{g}_{-1} \times \mathfrak{g}_{-1}^* \times \mathfrak{g}_0^*$, where $c$ and $d$ are skew-symmetric matrices with odd entries and $y$ is an arbitrary $n \times n$-matrix with even entries. For any $v \in V$ the image $r(v) \in \mathfrak{g}^*$ can be calculated as

$$(6.2)$$
$$\mathrm{Ad}^*_{\begin{pmatrix} 1_n & 0 \\ c & 1_n \end{pmatrix}} \begin{pmatrix} y & d \\ 0 & y^t \end{pmatrix} = \begin{pmatrix} 1_n & 0 \\ c & 1_n \end{pmatrix} \begin{pmatrix} y & d \\ 0 & y^t \end{pmatrix} \begin{pmatrix} 1_n & 0 \\ -c & 1_n \end{pmatrix} = \begin{pmatrix} y - dc & d \\ cy - y^tc - cdc & y^t + cd \end{pmatrix}$$

Denote the open subset $\mathfrak{g}_{-1} \times \mathfrak{g}_{-1}^* \times R_0 \subset \mathfrak{g}_{-1} \times \mathfrak{g}_{-1}^* \times \mathfrak{g}_0 \cong V$ by $R$.

**Lemma 6.3.** *(a)* $r : V \to \mathfrak{g}^*$ *is a $G$-equivariant map;*
*(b)* $\ker r^* = I_Q$;
*(c)* $r(R)$ *is dense in $Q$ and $r : R \to r(R)$ is an isomorphism.*

*Proof.* (a) is obvious from construction of $r$. Let us show (b). First, note that the underlying manifold of $V$ is $\mathfrak{g}_0^*$. As it is easy to check, the restriction of $r$ to $\mathfrak{g}_0^*$ is equal to the identity map $\mathfrak{g}_0^* \to \mathfrak{g}_0^*$. Let $\mu : G \times \mathfrak{g}_0^* \to V$ be the natural map induced by the action of $G$ on $V$ and the embedding $\mathfrak{g}_0^* \subset V$. We claim that the image of $\mu$ is dense in $V$. Indeed, it is sufficient to show that the image of $G_1 \times \mathfrak{g}_0^*$ is dense in the fiber $\mathfrak{g}_{-1}^* \times \mathfrak{g}_0^* \subset V$. But this follows from Lemma 6.1. Thus $\mu^*$ is injective. $G$-equivariance of $r$ implies that the map $\xi : G \times \mathfrak{g}_0^* \to \mathfrak{g}^*$ (defined in the begining of this section) equals $r \circ \mu$. Therefore by the injectivity of $\mu^*$, $I_Q = \operatorname{Ker} \xi^* = \operatorname{Ker} r^*$. To check $(c)$ note that $R$ is dense in $V$, hence $r(R)$ is dense in $Q$. To prove that $r : R \to r(R)$ is an isomorphism one calculates easily that $d_x r$ has a maximal rank for any regular point $x \in R_0$. □

By analogy with the even situation we consider $V$ as a resolution of singularities of $Q$. The mapping $r$ gives a "geometric description" of $Q$ as one can see from the following.

**Example 6.4.** Let $n = 2$. Denote the RHS of (6.2) by $\begin{pmatrix} q & s \\ u & q^t \end{pmatrix}$. Since $cdc = 0$ for $n = 2$, a straightforward calculation shows that $r(v)$ satisfies the equations

$$q_{12}u_{11} + q_{21}u_{22} = 2q_{21}u_{12} + (q_{11} - q_{22})\,u_{22} = 2q_{12}u_{12} + (q_{22} - q_{11})\,u_{11} = 0,$$
$$u_{11}u_{22} = u_{12}u_{11} = u_{12}u_{22} = 0.$$

One can check that the set of solutions of these equations coincides with $Q$. Note that in this case the equations do not involve $s$. Calculating the tangent spaces, we see that a point $x \in Q$ is singular iff $q$ is a scalar matrix.

We do not know equations on $Q$ for a general $n$.

Now let us formulate the following

**Conjecture 6.5.** The Poisson superalgebras $\operatorname{Gr} \bar{U}$ and $\mathcal{O}(Q)$ are isomorphic.

We give a proof of a weaker statement. Let $S \in U$ is defined as in Section 4, recall that its natural projection $\bar{S} \in \bar{U}$ belongs to the center of $\bar{U}$. Let $\widetilde{S}$ denote the the image of $S$ in $\operatorname{Gr} U$. Both algebras $\operatorname{Gr} \bar{U}$ and $\mathcal{O}(Q)$ are the quotients of $\operatorname{Gr} U$, and therefore the localized algebras $\operatorname{Gr} \bar{U} \left[ \widetilde{S}^{-1} \right]$ and $\mathcal{O}(Q) \left[ \widetilde{S}^{-1} \right]$ are both Poisson superalgebras.

**Theorem 6.6.** *The Poisson superalgebras* $\operatorname{Gr} \bar{U} \left[ \widetilde{S}^{-1} \right]$ *and* $\mathcal{O}(Q) \left[ \widetilde{S}^{-1} \right]$ *are isomorphic.*

*Proof.* Consider the following coinduced $U$-module

$$\mathcal{F} = \operatorname{Hom}_{U(\mathfrak{p})} (U, U(\mathfrak{g}_0)) .$$

Here the action of $U(\mathfrak{p})$ on $U(\mathfrak{g}_0)$ is defined by the conditions that $U(\mathfrak{g}_0)$ acts by left multiplication and $\mathfrak{g}_1$ acts by zero.

**Lemma 6.7.** $\operatorname{Ann}_U \mathcal{F} = I$.

*Proof.* Let $\mathcal{F}' = \mathrm{Hom}_{\bar{U}^+}\left(\bar{U}, U\left(\mathfrak{g}_0\right)\right)$. Since $\bar{U}$ is a quotient of $U$, there is the canonical embedding $\mathcal{F}' \subseteq \mathcal{F}$. We claim that $\mathcal{F}' = \mathcal{F}$. Indeed, by corollary 3.5 $U\left(\mathfrak{g}_0\right)$ is the subalgebra of $\bar{U}^+$ and $\bar{U} = \bar{U}^+ U\left(\mathfrak{g}_{-1}\right)$ is the decomposition with unique factorization. Recall also that $\mathfrak{g}_{-1}$ is supercommutative, hence $U\left(\mathfrak{g}_{-1}\right) = S\left(\mathfrak{g}_{-1}\right)$. Therefore we have an isomorphism of vector spaces.

$$(6.3) \qquad \mathcal{F}' \cong U\left(\mathfrak{g}_{-1}\right)^* \otimes U\left(\mathfrak{g}_0\right) = S\left(\mathfrak{g}_{-1}\right)^* \otimes U\left(\mathfrak{g}_0\right) \cong \mathcal{F}.$$

Hence $\mathcal{F}$ is well defined $\bar{U}$-module, thus $I \subseteq \mathrm{Ann}_U \mathcal{F}$.

It is easy to see that $S^d\left(\mathfrak{g}_{-1}\right)^* \otimes U\left(\mathfrak{g}_0\right) = \mathcal{F}^{\mathfrak{g}_1}$ and $\mathcal{F}^{\mathfrak{g}_1}$ generates $\mathcal{F}$ as $U$-module. Note that $\mathcal{F}^{\mathfrak{g}_1}$ is isomorphic to $U\left(\mathfrak{g}_0\right) \otimes C_{-\sigma}$ as $U\left(\mathfrak{g}_0\right)$-module, where $C_{-\sigma}$ is the one-dimensional $\mathfrak{g}_0$-module of weight $-\sigma$. Therefore $\mathcal{F} \cong \mathrm{Ind}\left(U\left(\mathfrak{g}_0\right) \otimes C_{-\sigma}\right)$. Hence any Verma module $M_\mu$ is a quotient of $\mathcal{F}$ and by Theorem 3.4 $\mathrm{Ann}_U \mathcal{F} \subseteq I$. $\qquad\square$

Let $\mathcal{D}$ be the associative algebra of differential operators on $K$ and $\bar{\mathcal{D}} = \mathcal{D} \otimes U\left(\mathfrak{g}_0\right)$. Since $\mathcal{O}\left(K\right) = S\left(\mathfrak{g}_{-1}^*\right)$ and $\mathcal{D}$ is a Clifford algebra, $\mathcal{D}$ is isomorphic to $\mathrm{End}_{\mathbb{C}}\left(S\left(\mathfrak{g}_{-1}^*\right)\right)$. Therefore $\bar{\mathcal{D}} \cong \mathrm{End}_{\mathbb{C}}\left(S\left(\mathfrak{g}_{-1}^*\right)\right) \otimes U\left(\mathfrak{g}_0\right)$. Using the identification $\mathcal{F} \cong S\left(\mathfrak{g}_{-1}\right)^* \otimes U\left(\mathfrak{g}_0\right)$ as in (6.3) define $\bar{\mathcal{D}}$ action on $\mathcal{F}$. Let $i \colon \bar{\mathcal{D}} \to \mathrm{End}\left(\mathcal{F}\right)$ be the homomorphism induced by this action. Obviously, $i$ is injective. Let $j : U \to \mathrm{End}\left(\mathcal{F}\right)$ be the homomorphism induced by the action of $U$ on $\mathcal{F}$. One can check that $j\left(U\right) \subseteq i\left(\bar{\mathcal{D}}\right)$, and this gives a homomorphism $\gamma = i^{-1} \circ j : U \to \bar{\mathcal{D}}$. Lemma 6.7 implies

**Corollary 6.8.** $\mathrm{Ker}\,\gamma = I$.

The algebra $\bar{\mathcal{D}}$ has the filtration induced by the natural filtration of $U\left(\mathfrak{g}_0\right)$ and the natural filtration of $\mathcal{D}$. We can identify $\mathrm{Gr}\,\bar{\mathcal{D}}$ with $\mathcal{O}\left(V\right)$. Indeed, the cotangent bundle $\mathcal{T}K^* \cong G \times_P \left(\mathfrak{g}/\mathfrak{p}\right)^*$ is a subbundle of $V$. The identification $V \cong \mathfrak{g}_{-1} \times \mathfrak{g}_{-1}^* \times \mathfrak{g}_0^*$ induces the isomorphisms $\mathcal{T}K^* \cong \mathfrak{g}_{-1} \times \mathfrak{g}_{-1}^*$ and $V \cong \mathcal{T}K^* \times \mathfrak{g}_0^*$. Therefore

$$\mathcal{O}\left(V\right) \cong \mathcal{O}\left(\mathcal{T}K^*\right) \otimes \mathcal{O}\left(\mathfrak{g}_0^*\right) \cong \mathrm{Gr}\,\mathcal{D} \otimes \mathrm{Gr}\,U\left(\mathfrak{g}_0\right) = \mathrm{Gr}\,\bar{\mathcal{D}}.$$

It is easy to check that $\gamma : U \to \bar{\mathcal{D}}$ is compatible with the filtrations on $U$ and $\bar{\mathcal{D}}$. Therefore $\gamma$ induces the homomorphism $\mathrm{Gr}\,\gamma : \mathrm{Gr}\,U \to \mathcal{O}\left(V\right)$ of graded algebras.

**Lemma 6.9.** $\mathrm{Gr}\,\gamma = r^*$.

*Proof.* The homomorphism $\mathrm{Gr}\,\gamma$ of supercommutative algebras induces a morphism $\theta : V \to \mathfrak{g}^*$ of affine supervarieties. We have to check that $\theta = r$. Consider the submanifold $V' = \mathfrak{g}_{-1}^* \times \mathfrak{g}_0^* \subset V$ (the fiber over $\left[\mathfrak{g}_{-1}\right] \in K$). By the definitions of $r$ and $\theta$, $r$ and $\theta$ coincide on $\mathcal{V}'$. Note that $\theta$ is $G$-equivariant as well as $r$. Since $GV' = V$, we obtain $\theta = r$. $\qquad\square$

Thus, $\mathrm{Gr}\,I \subset \mathrm{Ker}\,\mathrm{Gr}\,\gamma = I_Q$. Since $\bar{U} = U/\mathrm{Ker}\,\gamma$, $\gamma$ induces the injective homomorphism $\bar{\gamma} : \bar{U} \to \bar{\mathcal{D}}$, which is also compatible with the filtrations (here the filtration on $\bar{U}$ is the image filtration of $U$). Therefore $\bar{\gamma}$ induces the homomorphism $\mathrm{Gr}\,\bar{\gamma} \colon \mathrm{Gr}\,\bar{U} \to \mathcal{O}\left(V\right)$.

We will need again elements $T \in U$ and $S_0 \in U(\mathfrak{g}_0)$ defined in Section 4. Recall that for any $x \in U$ we denote by $\bar{x}$ its image in $\bar{U}$. Finally we deal with three filtered algebras $U, \bar{U}$ and $\bar{\mathcal{D}}$ and by $\deg x$ we denote the degree of an element $x$ of one of these algebras.

**Lemma 6.10.** (a) $\gamma(T) = \bar{\gamma}(\bar{T}) \in S^d(\mathfrak{g}_{-1}^*) \otimes S_0$;

(b) $\bar{\gamma}(\bar{U}\bar{T}\bar{U}) = \bar{\mathcal{D}}(1 \otimes S_0)$;

(c) $\bar{\gamma} : \bar{U}\bar{T}\bar{U} \to \bar{\mathcal{D}}(1 \otimes S_0)$ is an isomorphism preserving filtration.

*Proof.* Recall that $[z_0, \bar{T}] = d\bar{T}$ and $[\mathfrak{g}_0', \bar{T}] = 0$. Therefore $\gamma(T) \in S^d(\mathfrak{g}_{-1}^*) \otimes g_0$, where $g_0$ is some element in the center of $U(\mathfrak{g}_0)$. Furthermore $\deg T = 2d$, therefore $\deg \gamma(T) \leq 2d$, and $\deg g_0 \leq 2d$. On the other hand, $T$ acts by zero on any Verma module $M_\mu$ such that $\Delta(\mu) = 0$. Denote by $h(g_0)$ the image of $g_0$ under the Harish-Chandra map. Then $\Delta$ divides $h(g_0)$ and the degree of $h(g_0)$ is not higher than $2d$. This implies $h(g_0)$ is proportional to $\Theta$, i.e. (a) is proven.

To show (b) let $L_T$ be the ad $\mathfrak{g}$-submodule in $\bar{U}$ generated by $\bar{T}$. A straightforward calculation shows that $\bar{\gamma}(L_T) = S(\mathfrak{g}_{-1}^*) \otimes S_0$. Since $\mathfrak{g}_1\bar{T} = \bar{T}\mathfrak{g}_1 = 0$, we have $\bar{U}\bar{T}\bar{U} = U(\mathfrak{g}_0 \oplus \mathfrak{g}_{-1}) L_T$. Hence $\bar{\gamma}(\bar{U}\bar{T}\bar{U}) = \bar{\mathcal{D}}(1 \otimes S_0)$.

Let us prove (c). Injectivity of $\bar{\gamma}$ implies that $\bar{\gamma} : \bar{U}\bar{T}\bar{U} \to \bar{\mathcal{D}}(1 \otimes S_0)$ is an isomorphism, and it is left to check that $\bar{\gamma}$ preserves the filtrations. Let $\{x_i\}$ be a basis in $L_T$. For any $y \in \bar{U}\bar{T}\bar{U}$ one can find $g_i \in U(\mathfrak{g}_0 \oplus \mathfrak{g}_{-1})$ such that $y = \sum g_i x_i$. Obviously $\deg \gamma(g_i) = \deg g_i$ and $\deg \gamma(x_i) = \deg x_i = 2d$. Therefore

$$\deg \gamma(y) = \max\{\deg \gamma(g_i) \gamma(x_i)\} = \max\{\deg g_i\} + 2d \geq \deg y.$$

That proves (c). $\qquad\qquad\square$

**Corollary 6.11.** $\widetilde{S} \operatorname{Ker} \operatorname{Gr} \bar{\gamma} = 0$, and therefore the natural homomorphism $\operatorname{Gr} \bar{\gamma} : \operatorname{Gr} \bar{U}\left[\widetilde{S}^{-1}\right] \to \mathcal{O}(V)\left[S_0^{-1}\right]$ is an isomorphism.

*Proof.* By Lemma 6.10 (c) $\operatorname{Gr} \bar{\gamma} : \operatorname{Gr}(\bar{U}\bar{T}\bar{U}) \to \operatorname{Gr}(\bar{\mathcal{D}}(1 \otimes S_0))$ is an isomorphism of graded algebras. Therefore $\widetilde{S} \operatorname{Gr} \bar{U} \subset \operatorname{Gr}(\bar{S}\bar{U}) \subset \operatorname{Gr}(\bar{U}\bar{T}\bar{U})$ and thus $\operatorname{Ker} \operatorname{Gr} \bar{\gamma} \cap \widetilde{S} \operatorname{Gr} \bar{U} = \{0\}$. That implies $\widetilde{S} \operatorname{Ker} \operatorname{Gr} \bar{\gamma} = 0$. $\qquad\square$

As follows from Lemma 6.9 $\operatorname{Gr} \bar{\gamma} = \bar{r}^*$, where $\bar{r}^* : \mathcal{O}(Q) \to \mathcal{O}(V)$ is induced by $r^*$. Therefore Corollary 6.11 implies the theorem. $\qquad\square$

*Remark* 6.12. To prove the conjecture we have to show that $\widetilde{S}$ is not a zero divisor in $\operatorname{Gr} \bar{U}$, which means that the supervariety $\mathcal{X} \subset \mathfrak{g}^*$ associated to $\operatorname{Gr} I$ is "irreducible". Theorem 6.6 in this terms means that $Q$ is an "irreducible component" of $\mathcal{X}$. It also implies that $\mathcal{X}$ and $Q$ coincide after removing singularities.

Finally let us note that Lemma 6.10 (c) implies the following statement.

**Corollary 6.13.** The localized algebras $\bar{U}\left[\bar{S}^{-1}\right]$ and $\operatorname{End}_{\mathbb{C}} U(\mathfrak{g}_{-1}) \otimes U(\mathfrak{g}_0)\left[S_0^{-1}\right]$ are isomorphic.

## References

1. M. Gorelik, *The center of a simple p-type Lie superalgebra*, Journal of Algebra.
2. A. Joseph, *Kostant's problem, Goldie rank and the Gelfand-Kirillov conjecture*, Inventiones Mathematicae **56** (1980), no. 3, 191–213.
3. V. G. Kac, *Lie superalgebras*, Adv. Math. **26** (1977), 8–96.
4. E. Kirkman and J. Kuzmanovich, *Minimal prime ideals in enveloping algebras of Lie superalgebras*, Proceedings of American Mathematical Society **124** (1996), 1693–1702.
5. Yu. Manin, *Gauge field theory and complex geometry*, Grundlehren der Mathematischen Wissenschaften [Fundamental Principles of Mathematical Sciences], vol. 289, Springer-Verlag, Berlin, 1988, Translated from the Russian by N. Koblitz and J. R. King.
6. I. Penkov, *Borel-Weil-Bott theory for classical Lie supergroups*, J. Soviet Math.**51**,(1988),2108–2140.

DEPT. OF MATHEMATICS, UNIVERSITY OF CALIFORNIA AT BERKELEY, BERKELEY, CA 94720
*E-mail address*: sergeant@math.berkeley.eau